# SADDLE-POINT THEOREMS FOR THE GENERALIZED CONE-CONVEX OPTIMIZATIONS OF SET-VALUED FUNCTIONS

Renying Zeng

School of Mathematical Sciences, Chongqing Normal University, Chongqing, China
Email: renying.zeng@saskpolytech.ca

**Abstract**: This paper works with preconvexlike set-valued vector optimization problems in topological linear spaces. A Fakas-Minkowski alternative theorem, a scalarization theorem, some vector saddle-point theorems and some scalar saddle point theorem are proved.

## 1. Introduction and Preliminary Results

Convexity and generalized convexities have played very important roles in optimization (including mathematical programming, and multi-objective mathematical programming). Convexity and generalized convexities are also very important in other areas such as multi-criteria decision, non-smooth analysis, and control theory, etc. In this paper, we work with preconvexlike vector optimization problems of set-valued mappings in topological linear spaces.

Let $X$ be a real topological linear space, $X^*$ the topology dual of $X$. A subset $X_+$ of $X$ is said to be a convex cone if

$$\alpha x^1 + \beta x^2 \in X_+, \forall x^1, x^2 \in X_+, \ \forall \alpha, \beta \geq 0 .$$

A real topological linear space $Y$ with a convex cone is said to be an ordered topological linera space. We denote $intX_+$ the topological interior of $X_+$. The partial order on $X$ is defined by

$$x^1 \leq x^2, \text{ if } x^1 - x^2 \in X_+,$$
$$x^1 < x^2, \text{ if } x^1 - x^2 \in int\, X_+.$$

The subset of $X^*$

$$X_+^* = \{\xi \in X^* : \langle x, \xi \rangle \geq 0, \forall x \in X_+\}$$

is said to be the dual cone of the cone $X_+$, where $\langle x, \xi \rangle = \xi(x)$.

Suppose that $X$ and $Y$ are two real topological liner spaces. Let $f$: $X \to 2^Y$ be a set-valued function, where $2^Y$ denotes the power set of $Y$.

Let $D$ be a nonempty subset of $X$. Setting

$$f(D) = \bigcup_{x \in D} f(x),$$
$$\langle f(x), \eta \rangle = \{\langle y, \eta \rangle : y \in f(x)\},$$
$$\langle f(D), \eta \rangle = \bigcup_{x \in D} \langle f(x), \eta \rangle.$$

For $x \in D, \eta \in Y^*$, write

$$\langle f(x), \eta \rangle \geq 0 \text{, if } \langle y, \eta \rangle \geq 0, \forall y \in f(x),$$
$$\langle f(D), \eta \rangle \geq 0 \text{, if } \langle f(x), \eta \rangle \geq 0, \forall x \in D.$$

A set-valued function $f$: $X \to 2^Y$ is said to be $Y_+$ -convex on $D$ if $\forall x^1, x^2 \in D$, $\forall \alpha \in [0, 1]$, one has

$$\alpha f(x^1) + (1-\alpha) f(x^2) \subseteq f(\alpha x^1 + (1-\alpha) x^2) + Y_+.$$

A set-valued function $f$: $X \to 2^Y$ is said to be $Y_+$ -convexlike on $D$ if $\forall x^1, x^2 \in D$, $\forall \alpha \in [0, 1]$, $\exists x^3 \in D$ such that

$$\alpha f(x^1) + (1-\alpha) f(x^2) \subseteq f(x^3) + Y_+.$$

We introduce below the concept of preconvexlike about vector-valued functions to set-valued functions.

**Definition 1.1** A set-valued function $f: X \to 2^Y$ is said to be $Y_+$ - preconvexlike on $D$ if $\forall x^1, x^2 \in D$, $\forall \alpha \in (0, 1)$, $\exists x^3 \in D$, and $\exists \tau > 0$, such that

$$\alpha f(x^1) + (1-\alpha) f(x^2) \subseteq \tau f(x^3) + Y_+.$$

Any set- and real-valued function is convexlike, it is known that a set-valued convexlike function is not necessary to be a convex function.

Let $X = Y = R^2$ (the two dimensional Euclidean space), $Y_+$ be the first quadrant of $R^2$. Example 1.1 illustrates that a set-valued preconvexlike function is not necessary to be a convexlike function.

**Example 1.1** Take $D = [Y_+ \setminus \{(x_1, x_2) : 0 \leq x_1 \leq 1, 0 \leq x_2 \leq 1\}] \cup \{(0,1), (1,0)\}$, and define a set-valued function $f: X \to 2^Y$ by

$$f(x_1, x_2) = D.$$

Then $f$ is preconvexlike not convexlike.

We also introduce the concept of preaffine set-valued functions here, which extends the definition of affine functions.

**Definition 1.2** A set-valued function $f: X \rightarrow 2^Y$ is said to be preaffine on $D$ if $\forall x^1$, $x^2 \in D$, $\forall \alpha \in (0, 1)$, $\exists x^3 \in D$, and $\exists \tau > 0$, such that

$$\alpha f(x^1) + (1-\alpha) f(x^2) \subseteq \tau f(x^3).$$

Assume that $X$, $Y$, $Z$ and $W$ are topological linear spaces with convex cones $X_+$, $Y_+$, $Z_+$, and $W_+$ respectively; $X_+$, $Y_+$, $Z_+$, and $W_+$ have nonempty interiors *int* $X_+$, *int* $Y_+$, *int* $Z_+$, and *int* $W_+$.

**Lemma 1.1** [12] Suppose that $Y$ is a topological linear space. If $\xi \in X_+^* \setminus \{O\}$, and if $x^0 \in int\, X_+$, then $\langle x^0, \xi \rangle > 0$.

**Theorem 1.2** (Fakas-Minkowski Theorem of Alternate) Suppose

(a) $f: X \rightarrow 2^Y$ is $Y_+$ - preconvexlike on $D$, and $g: Y \rightarrow 2^Z$ is $Z_+$ - preconvexlike on $D$,

(b) $h: X \rightarrow 2^W$ is preaffine on $D$, and $int[\, h(D)] \neq \emptyset$,

and, (i) and (ii) denote the systems

(i) $\exists x^0 \in D$ such that $f(x^0) \cap (-\text{int}\, Y_+) \neq \emptyset, g(x^0) \cap (-Z_+) \neq \emptyset$, and $O \in h(x^0)$, where $O$ is the zero element of a topological linear space;

(ii) $\exists (\xi, \eta, \zeta) \in Y_+^* \times Z_+^* \times W^*$, $(\xi, \eta, \zeta) \neq O$, such that

$$\xi(f(x)) + \eta(g(x)) + \varsigma(h(x)) \geq 0, \forall x \in D.$$

If (i) has no solution then (ii) has a solution.

If (ii) has a solution with $\xi \neq O$, then (i) has no solution.

**Proof.** At first, aim to show that the Cartesian product

$$A = [\cup_{t>0} tf(D) + \text{int}\, Y_+] \times [\cup_{t>0} tg(D) + Y_+] \times [\cup_{t>0} th(D)]$$

is convex.

Take $$c^i = (t_i' f(x^i) + y^i, t_i'' g(x^i) + z^i, t_i''' h(x^i)) \in A$$

where $(y^i, z^i) \in (intY_+) \times Z_+$, $t_i', t_i'', t_i'''$ are positive numbers ($i$ = 1, 2,3). For $\alpha \in (0, 1)$, set $y^0 = \alpha y^1 + (1-\alpha) y^2$, $z^0 = \alpha z^1 + (1-\alpha) z^2$. Noting that $intY_+$ and $Z_+$ are convex cones, we have $y^0 \in intY_+$, $z^0 \in Z_+$. Since $f$ and $g$ are preconvexlike and $g$ is preaffine, $\exists\, x_3', x_3'', x_3''' \in D$, and $\exists\, \tau', \tau'', \tau''' > 0$, such that

$$\frac{\alpha t_1'}{\alpha t_1'+(1-\alpha)t_2'}f(x^1)+\frac{(1-\alpha)t_2'}{\alpha t_1'+(1-\alpha)t_2'}f(x^2)-\tau' f(x_3'):=y'\in Y_+,$$
$$\frac{\alpha t_1'}{\alpha t_1''+(1-\alpha)t_2''}g(x^1)+\frac{(1-\alpha)t_2''}{\alpha t_1''+(1-\alpha)t_2''}g(x^2)-\tau'' g(x_3''):=z'\in Z_+,$$
$$\frac{\alpha t_1'''}{\alpha t_1'''+(1-\alpha)t_2'''}h(x^1)+\frac{(1-\alpha)t_2'''}{\alpha t_1'''+(1-\alpha)t_2'''}h(x^2)-\tau''' h(x_3''')=O$$

But

$$\begin{aligned}
&\alpha t_1' f(x^1)+(1-\alpha)t_2' f(x^2)+\alpha y^1+(1-\alpha)y^2\\
&=(\alpha t_1'+(1-\alpha)t_2')[\frac{\alpha t_1}{\alpha t_1'+(1-\alpha)t_2'}f(x^1)+\frac{(1-\alpha)t_2'}{\alpha t_1'+(1-\alpha)t_2'}f(x^2)]+y^0\\
&=(\alpha t_1'+(1-\alpha)t_2')[y'+\tau' f(x_3')]+y^0 \qquad (1.1)\\
&=(\alpha t_1'+(1-\alpha)t_2')\tau' f(x_3')+(\alpha t_1'+(1-\alpha)t_2')y'+y^0\\
&\subseteq t' f(x_3')+Y_+ +\operatorname{int}Y_+\subseteq t' f(x_3')+\operatorname{int}Y_+.
\end{aligned}$$

Similarly

$$\begin{aligned}
&\alpha t_1'' g(x^1)+(1-\alpha)t_2'' g(x^2)+\alpha z^1+(1-\alpha)z^2\\
&\subseteq t'' g(x_3'')+Z_+ +Z_+\subseteq t'' g(x_3'')+Z_+.
\end{aligned} \qquad (1.2)$$

And

$$\alpha t_1''' h(x^1)+(1-\alpha)t_2''' h(x^2)\subseteq t''' h(x_3''')\subseteq t''' h(x_3'''). \qquad (1.3)$$

From (1.1), (1.2) and (1.3), if $c^1, c^2\in A$, we have

$$\alpha c^1+(1-\alpha)c^2\in[t' f(x_3')+\operatorname{int}Y_+]\times[t'' g(x_3'')+Z_+]\times[t''' h(x_3''')]\subseteq A$$

Which means $A$ is convex.

Secondly, suppose that (i) has no solution, and want to show (ii) has a solution.

Because (i) has no solution, $O\notin A$.

And because we assume that $intY_+\neq\varnothing$, $intZ_+\neq\varnothing$, and $int[h(D)]\neq\varnothing$, we have $int\, A\neq\varnothing$.

By the separation theorem of convex sets of a topological linear space

$$\exists\ \gamma=(\xi,\eta,\zeta)\in Y^*\times Z^*\times W^*,\ \gamma=(\xi,\eta,\zeta)\neq O, \qquad (1.4)$$

such that

$$\gamma(A)\geq 0,$$

i.e., $\forall t', t'', t'''>0, \forall y^0\in\operatorname{int}Y_+, \forall z^0\in Z_+$, we have

$$\left\langle t' f(x) + y^0, \xi \right\rangle + \left\langle t'' g(x) + z^0, \eta \right\rangle + \left\langle t''' h(x), \varsigma \right\rangle \geq 0, \forall x \in D.$$

Noting that $\forall y \in Y_+, \exists y_n \in \operatorname{int} Y_+, (n = 1,2....)$ such that $y_n \to y (n \to \infty)$, and $\forall s' > 0, s'' > 0, s' y_n \in \operatorname{int} Y_+, s'' z^0 \in Z_+$, we have

$$\left\langle t' f(x) + s' y_n, \xi \right\rangle + \left\langle t'' g(x) + s'' z^0, \eta \right\rangle + \left\langle t''' h(x), \varsigma \right\rangle \geq 0, \forall x \in D. \qquad (1.5)$$

Let $t' = t'' = s' = s'' = 1, t''' = 0, n \to \infty$, we get $\left\langle y, \xi \right\rangle + \left\langle z^0, \eta \right\rangle \geq 0, \forall y \in Y_+, \forall z^0 \in Z_+$. Therefore

$$(\xi, \eta, \zeta) \in Y_+^* \times Z_+^* \times W^*.$$

Let $t' = t'' = t''' = 1$, $s', s'' \to 0$ in (1.5) we will get

$$\left\langle f(x), \xi \right\rangle + \left\langle g(x), \eta \right\rangle + \left\langle h(x), \varsigma \right\rangle \geq 0, \forall x \in D.$$

Which shows us that (ii) has a solution.

Thirdly, assume that (ii) has a solution $(\xi, \eta, \zeta)$ with $\xi \neq O$, i.e.

$$\left\langle f(x), \xi \right\rangle + \left\langle g(x), \eta \right\rangle + \left\langle h(x), \varsigma \right\rangle \geq 0, \forall x \in D.$$

We are going to prove (i) has no solution.

Otherwise, if (i) has a solution $x \in D$, there would exist $y \in f(x), z \in g(x)$, and $w \in h(x)$ such that $y \in -\operatorname{int} Y_+, z \in -Z_+, w = O$. By Lemma 1.1 we would have

$$\left\langle y, \xi \right\rangle + \left\langle z, \eta \right\rangle + \left\langle w, \varsigma \right\rangle < 0, \forall x \in D.$$

Which is a contradiction. □

## 2. Scalarization

Consider the following vector optimization problem with set-valued functions:

$$\text{(VP)} \qquad \begin{aligned} & Y_+ - \min f(x), \\ & s.t., g(x) \cap (-Z_+) \neq \varnothing, 0_Z \in h(x), \\ & x \in X. \end{aligned}$$

Let $D$ be the feasible set of (VP), i.e., $D = \{x \in X : g(x) \cap (-Y_+) \neq \varnothing, 0_W \in h(x)\}$.

**Definition 2.1** $\bar{x} \in D$ is said to be a weakly efficient solution of (VP) if $\exists \bar{y} \in f(\bar{x})$ such that $(\bar{y} - f(D)) \cap \operatorname{int} Y_+ = \varnothing$.

**Definition 2.2** The problem (VP) is said to satisfy the Slater constraint qualification (SC) if $\forall (\eta, \varsigma) \in (Z_+^* \times W^*) \setminus \{O\}$, $\exists x \in D$ such that

$$(-R_+) \bigcap [\eta(g(x)) \cap \varsigma(h(x))] \neq \varnothing,$$

where $R_+$ is the set of all positive real numbers.

**Definition 2.3** $\bar{x} \in D$ is said to be an optimal solution of the scalar optimization problem (VPS), if $\exists \bar{y} \in f(\bar{x})$, and $\exists\, y_+ \in Y_+$ such that

$$\xi(f(x)) \subseteq \xi(\bar{y}) + Y_+, \forall x \in D.$$

By use of the following Theorem of Scalarization, we may convert a vector optimization problem into a scalar optimization problem.

**Theorem 2.1** (Scalarization) Suppose $\bar{x} \in D$, and

(a) $f(x) - f(\bar{x}) : X \to 2^Y$ is $Y_+$-preconvexlike on $D$, and $g(x)$: $Y \to 2^Z$ is $Z_+$-preconvexlike on $D$,

(b) $h(x)$: $X \to 2^W$ is preaffine on $D$, and

(c) (VP) satisfies the Slater constraint qualification (SC),

then $\bar{x}$ is a weakly efficient solution of (VP) if and only if $\exists \xi \in Y_+^* \setminus \{O\}$ such that $\bar{x}$ is an optimal solution of the following scalar optimization problem (VPS):

$$\text{(VPS)} \qquad \min_{x \in D} \xi(f(x)).$$

**Proof.** If $\exists \xi \in Y_+^* \setminus \{O\}$ and if $\bar{x} \in D$ is an optimal solution of the scalar optimization problem (VPS), then $\exists \bar{y} \in f(\bar{x})$, such that $\xi(f(x)) \subseteq \xi(\bar{y}) + Y_+, \forall x \in D$**.** So

$$(\bar{y} - f(D)) \cap \operatorname{int} Y_+ = \varnothing.$$

Therefore $\bar{x}$ is a weakly efficient solution of (VP).

On the other hand, suppose $\bar{x}$ is a weakly efficient solution of (VP), we want to show that $\bar{x}$ is an optimal solution of the scalar optimization problem (VPS).

From Definition 2.1 $\exists \bar{y} \in f(\bar{x})$ such that the following system

$$(f(x) - \bar{y}) \cap (-\operatorname{int} Y_+) \neq \varnothing, g(x) \cap (-Z_+) \neq \varnothing, O \in h(x)$$

has no solution for $x \in D$. Hence, Theorem 1.2 implies

$$\exists \xi \in Y_+^*, \eta \in Z_+^*, \varsigma \in W^* \text{ with } (\xi, \eta, \varsigma) \neq O \qquad (2.1)$$

such that

$$\xi(f(x) - \bar{y}) + \eta(g(x)) + \varsigma(h(x)) \geq 0, \forall x \in D.$$

i.e.,

$$\xi(f(x)) + \eta(g(x)) + \varsigma(h(x)) \geq \xi(\bar{y}), \forall x \in D. \qquad (2.2)$$

If $\xi = O$, then by (2.1) we get $(\xi, \eta) \neq O$. And (2.2) yields

$$\eta(g(x)) + \varsigma(h(x)) \geq 0, \forall x \in D.$$

This is contradicting to the Slater constraint qualification (SC). Therefore $\xi \neq O$. Hence, from $\bar{x} \in D$, i.e., $g(x) \cap (-Z_+) \neq \varnothing, O \in h(x)$, (2.1) gives

$$\xi(f(x)) \geq \xi(\bar{y}), \forall x \in D.$$

Which means $\bar{x}$ is an optimal solution of (VPS). □

## 3. Vector Saddle-Point Theorems

Write

$$P\min[A, Y_+] = \{y \in A : (y - A) \cap \operatorname{int} Y_+ = \varnothing\},$$
$$P\max[A, Y_+] = \{y \in A : (A - y) \cap \operatorname{int} Y_+ = \varnothing\}.$$

We familiar with that

**Lemma 3.1** $\bar{x} \in D$ is a weakly efficient solution of (VP), if and only if

$$f(\bar{x}) \cap P\min[f(D), Y_+] \neq \varnothing.$$

**Definition 3.1** A triple $(\bar{x}, \bar{S}, \bar{T}) \in X \times B^+(Z, Y) \times B(W, Y)$ is said to be a vector saddle-point of L if

$$L(\bar{x}, \bar{S}, \bar{T}) \cap P\min[L(X, \bar{S}, \bar{T}), Y_+] \cap P\max[L(\bar{x}, B^+(Z, Y), B(W, Y)), Y_+] \neq \varnothing.$$

Where

$$L(\bar{x}, \bar{S}, \bar{T}) = f(\bar{x}) + \bar{S}(g(\bar{x})) + \bar{T}(h(\bar{x})).$$

**Definition 3.1** A convex cone $Y_+$ is said to be pointed if $Y_+ \bigcap (-Y_+) = \{O\}$.

From now on, suppose all convex cones are pointed and closed.

**Theorem 3.1** A triple $(\bar{x},\bar{S},\bar{T}) \in X \times B^{+}(Z,Y) \times B(W,Y)$ is a vector saddle-point of L, if and only if $\exists \bar{y} \in f(\bar{x}), \bar{z} \in g(\bar{x})$, such that

(i) $\bar{y} \in P\min[L(X,\bar{S},\bar{T}), Y_{+}]$,
(ii) $g(\bar{x}) \subset -Z_{+}, h(\bar{x}) = \{O\}$,
(iii) $(f(\bar{x}) - \bar{y} - \bar{S}(\bar{z})) \cap \operatorname{int} Y_{+} = \varnothing$.

**Proof**. The necessity. Assume that $(\bar{x},\bar{S},\bar{T}) \in X \times B^{+}(Z,Y) \times B(W,Y)$ is a vector saddle-point of $L$. From Definition 3.1

$$L(\bar{x},\bar{S},\bar{T}) \cap P\min[L(X,\bar{S},\bar{T}), Y_{+}] \cap P\max[L(\bar{x}, B^{+}(Z,Y) \times B(W,Y)), Y_{+}] \neq \varnothing.$$

So, $\exists \bar{y} \in f(\bar{x}), \bar{z} \in g(\bar{x}), \bar{w} \in h(\bar{x})$, i.e.,

$$\bar{y} + S(\bar{z}) + T(\bar{w}) \in L(\bar{x},\bar{S},\bar{T}) = f(\bar{x}) + S(g(\bar{x})) + T(h(\bar{x})), \tag{3.1}$$

such that

$$\begin{aligned} &\{ f(\bar{x}) + S(g(\bar{x})) + T(h(\bar{x})) - [\bar{y} + \bar{S}(\bar{z}) + \bar{T}(\bar{w})] \} \cap \operatorname{int} Y_{+} = \varnothing, \\ &\forall (S,T) \in B^{+}(Z,Y) \times B(W,Y), \end{aligned} \tag{3.2}$$

and

$$(\bar{y} + \bar{S}(\bar{z}) + \bar{T}(\bar{w}) - [f(X) + \bar{S}(g(X)) + \bar{T}(h(X))]) \cap \operatorname{int} Y_{+} = \varnothing. \tag{3.3}$$

Taking $T = \bar{T}$ in (3.2) we get

$$S(z) - \bar{S}(\bar{z}) \notin \operatorname{int} Y_{+}, \forall z \in g(\bar{x}), \forall S \in B^{+}(Z,Y). \tag{3.4}$$

Aim to show that $-\bar{z} \in Z_{+}$.
Otherwise, since $O \in -Z_{+}$, if $-\bar{z} \notin Z_{+}$, we would have $-\bar{z} \neq O$,
Because $Z_{+}$ is a closed convex set, by the separate theorem $\exists \eta \in Z^{*} \setminus \{O\}$

$$\eta(tz_{+}) > \eta(-\bar{z}), \forall z \in Z_{+}, \forall t > 0. \tag{3.5}$$

i.e.,

$$\eta(z_{+}) > \frac{1}{t}\eta(-\bar{z}), \forall z \in Z_{+}, \forall t > 0.$$

Let $t \to \infty$ we obtain $\eta(z_{+}) \geq 0, \forall z \in Z_{+}$. Which means that $\eta \in Z_{+}^{*} \setminus \{O\}$. Meanwhile, $O \in Z_{+}$ and (3.5) yield that $\eta(\bar{z}) > 0$. Given $\tilde{z} \in \operatorname{int} Z_{+}$ and let

$$S(z)=\frac{\eta(z)}{\eta(\bar{z})}\tilde{z}+\bar{S}(z).$$

Then $\bar{S}\in B^{+}(Z,Y)$ and

$$S(\bar{z})-\bar{S}(\bar{z})=\tilde{z}\in \operatorname{int}Y_{+}.$$

Contradicting to (3.4). Therefore

$$-\bar{z}\in Z_{+}.$$

Now, aim to prove that $-g(\bar{x})\subseteq Z_{+}$.

Otherwise, if $-g(\bar{x})\not\subset Z_{+}$, then $\exists z_0\in g(\bar{x})$ such that $O\neq -z_0\notin Z_{+}$. Similar to the above $\exists\,\eta_0\in Z^{*}\setminus\{O\}$ such that $\eta_0\in Z_{+}^{*}\setminus\{O\}$, $\eta_0(z_0)>0$. Given $\tilde{z}\in \operatorname{int}Z_{+}$ and let

$$S_0(z)=\frac{\eta_0(z)}{\eta_0(z_0)}\tilde{z}.$$

Then $S_0\in B^{+}(Z,Y)$ and $S_0(z_0)=\tilde{z}\in \operatorname{int}Y_{+}$. And we have proved that $-\bar{z}\in Z_{+}$, so $-\bar{S}(\bar{z})\in Y_{+}$. Therefore

$$S_0(z_0)-\bar{S}(\bar{z})\in \operatorname{int}Y_{+}+Y_{+}\subseteq \operatorname{int}Y_{+}.$$

Again, contradicting to (3.4).

Therefore $-g(\bar{x})\subseteq Z_{+}$. Similarly, one has $-h(\bar{x})\subseteq W_{+}$. From (3.2) we get

$$[T(h(\bar{x}))-\bar{T}(\bar{w})]\cap \operatorname{int}Y_{+}=\varnothing.$$

Hence

$$T(\bar{w})-\bar{T}(\bar{w})\notin \operatorname{int}Y_{+},\forall T\in B(W,Y). \tag{3.6}$$

Similarly, from (3.2) again we have

$$T(w)-\bar{T}(\bar{w})\notin \operatorname{int}Y_{+},\forall w\in h(\bar{x}),\forall T\in B(W,Y). \tag{3.7}$$

If $\bar{w}\neq O$, since $-h(\bar{x})\subseteq W_{+}$ and $W_{+}$ is a pointed cone, we have $\bar{w}\notin W_{+}$. Because $Y_{+}$ is a closed convex set, by the separation theorem $\exists\varsigma\in W^{*}$, such that

$$\varsigma(w)<\varsigma(\bar{w}),\forall w\in W_{+}. \tag{3.8}$$

So $\varsigma(\bar{w})\neq 0$ since $O\in W_{+}$. Taking $y^{0}\in \operatorname{int}Y_{+}$ and define $T^{0}\in B^{+}(W,Y)$ by

$$T^0(w) = \frac{\varsigma(w)}{\varsigma(\overline{w})} y^0 + \overline{T}(w).$$

Then

$$T^0(\overline{w}) - \overline{T}(\overline{w}) = y^0 \in \operatorname{int} Y_+,$$

Contradicting to (3.6). Therefore $\overline{w} = O$.Thus

$$O \in h(\overline{x}).$$

Now, we'd like to prove $h(\overline{x}) = \{O\}$.

Otherwise, if $w^0 \in h(\overline{x}) : w^0 \neq O$, similar to (3.8) $\exists \varsigma^0 \in W^*$, such that $\varsigma^0(w) < \varsigma^0(w^0), \forall w \in W_+$. So $\varsigma^0(w^0) \neq 0$. Given $y_0 \in \operatorname{int} Y_+$ and define $T_0 \in B(W,Y)$, by

$$T_0(w) = \frac{\varsigma^0(w)}{\varsigma^0(w^0)} y_0.$$

Then $T_0(w^0) = y_0 \in \operatorname{int} Y_+$, i.e., $\overline{w} = O$ ) $T^0(w^0) - \overline{T}(\overline{w}) \in \operatorname{int} Y_+$. Contradicting to (3.7). Therefore we must have

$$h(\overline{x}) = \{O\}. \tag{3.9}$$

Combining (3.2), (3.3), (3.9), and we conclude that

$$\overline{y} \in P\min[L(X, \overline{S}, \overline{T}), Y_+], \tag{3.10}$$

and

$$(f(\overline{x}) - \overline{y} - \overline{S}(\overline{z})) \cap \operatorname{int} Y_+ = \varnothing.$$

We have proved that, if $(\overline{x}, \overline{S}, \overline{T}) \in X \times B^+(Z,Y) \times B(W,Y)$ is a vector saddle-point of $L$, then the conditions (i)-(iii) hold.

The sufficiency. Suppose that the conditions (i)-(iii) are satisfied. Note that $-g(\overline{x}) \subseteq Z_+, h(\overline{x}) = \{O\}$ means

$$-S(g(\overline{x})) \subseteq Y_+,\ T(h(\overline{x})) = \{O\},\ \forall (S,T) \in B^+(Z,Y) \times B(W,Y), \tag{3.11}$$

and the condition (i) states that

$$\{\overline{y} - [f(X) + \overline{S}(g(X)) + \overline{T}(h(X))]\} \cap \operatorname{int} Y_+ = \varnothing,$$

So $Y_+ + \text{int}\, Y_+ \subseteq Y_+$ and $-S(\bar{z}) \in Y_+$ together imply

$$\{\bar{y} + \bar{S}(\bar{z}) + \bar{T}(\bar{w}) - [f(X) + \bar{S}(g(X)) + \bar{T}(h(X))]\} \cap \text{int}\, Y_+ = \varnothing .$$

Hence

$$\bar{y} + \bar{S}(\bar{z}) + \bar{T}(\bar{w}) \in P\min[L(X, \bar{S}, \bar{T}), Y_+] .$$

On the other hand, since $(f(\bar{x}) - [\bar{y} + S(\bar{z})]) \cap \text{int}\, Y_+ = \varnothing$, from (3.11), and from $\text{int}\, Y_+ + Y_+ \subseteq \text{int}\, Y_+$ we conclude that

$$\{\bigcup_{(S,T)\in B^+(Z,Y)\times B(W,Y)} [f(\bar{x}) + S(g(\bar{x})) + T(h(\bar{x}))] - [\bar{y} + \bar{S}(\bar{z}) + \bar{T}(\bar{w})]\} \cap \text{int}\, Y_+ = \varnothing .$$

Hence

$$\bar{y} + \bar{S}(\bar{z}) + \bar{T}(\bar{w}) \in P\max[L(\bar{x}, B^+(Z,Y), B(W,Y)), Y_+]$$

Consequently,

$$L(\bar{x}, \bar{S}, \bar{T}) \cap P\min[L(X, \bar{S}, \bar{T}), Y_+] \cap P\max[L(\bar{x}, B^+(Z,Y), B(W,Y)), Y_+] \neq \varnothing.$$

Therefore $(\bar{x}, \bar{S}, \bar{T}) \in X \times B^+(Z,Y) \times B(W,Y)$ is a vector saddle-point of *L*. □

**Theorem 3.2** If $(\bar{x}, \bar{S}, \bar{T}) \in X \times B^+(Z,Y) \times B(W,Y)$ is a vector saddle-point of *L*, and if $O \in \bar{S}(g(\bar{x}))$, then $\bar{x}$ is a weak efficient solution of (VP).

**Proof**. Assume that $(\bar{x}, \bar{S}, \bar{T}) \in D \times B^+(Z,Y) \times B(W,Y)$ is a vector saddle-point of *L*, from Theorem 3.1 we have

$$-S(g(\bar{x})) \subseteq Y_+ , h(\bar{x}) = \{O\} . \qquad (3.12)$$

So $\bar{x} \in D$ (the feasible solution of (VP). And $\exists \bar{y} \in f(\bar{x})$ such that $\bar{y} \in P\min[L(X, \bar{S}, \bar{T}), Y_+]$, i.e.

$$(\bar{y} - [f(X) + \bar{S}(g(X)) + \bar{T}(h(X)]) \cap \text{int}\, Y_+ = \varnothing .$$

Thus

$$(\bar{y} - [f(D) + \bar{S}(g(\bar{x})) + \bar{T}(h(\bar{x})]) \cap \text{int}\, Y_+ = \varnothing .$$

By (3.12) and note that $O \in \bar{S}(g(\bar{x}))$ we get

$$(\bar{y} - f(D)) \cap \operatorname{int} Y_+ = \varnothing .$$

Therefore by Definition 2.1 $\bar{x}$ is a weakly efficient solution of (VP). □

**Theorem 3.3** Suppose $\bar{x} \in D$ is a weakly efficient solution of (VP), and

(a) *f*: $X \to 2^Y$ is $Y_+$ - preconvexlike on *D*, and *g*: $Y \to 2^Z$ is $Z_+$ - preconvexlike on *D*,
(b) *h*: $X \to 2^W$ is preaffine on *D*, and $int[\, h(D)] \neq \varnothing$,
(c) (VP) satisfies the Slater constrained qualification (SC).

If $-g(\bar{x}) \subseteq Z_+$, $h(\bar{x}) = \{O\}$, and if $\exists \bar{y} \in f(\bar{x})$ for which $(f(\bar{x}) - \bar{y}) \cap \operatorname{int} Y_+ = \varnothing$, then $\exists (\bar{S}, \bar{T}) \in B^+(Z,Y) \times B(W,Y)$ such that $(\bar{x}, \bar{S}, \bar{T}) \in X \times B^+(Z,Y) \times B(W,Y)$ is a vector saddle-point of L and $O \in \bar{S}(g(\bar{x}))$.

**Proof**. Assume that $\bar{x} \in D$ is a weakly efficient solution of (VP), then $\exists \bar{y} \in f(\bar{x})$ for which there is not any $x \in D$ such that $f(x) - \bar{y} \in -\operatorname{int} Y_+$. So, there is not any $x \in X$ such that

$$f(x) - \bar{y} \in -\operatorname{int} Y_+, g(x) \in -Z_+, O \in h(x).$$

By Theorem 1.2 $\exists (\xi, \eta, \varsigma) \in Y_+^* \times Z_+^* \times W^* \setminus \{O\}$ such that

$$\xi(f(x) - \bar{y}) + \eta(g(x)) + \varsigma(h(x)) \geq 0, \forall x \in D. \qquad (3.13)$$

Since $\bar{y} \in f(\bar{x})$ and $O \in h(\bar{x})$, take $x = \bar{x}$ in the above we obtain $\eta(g(\bar{x})) \geq 0$. And $\bar{x} \in D$ and $\eta \in Z_+^*$ imply that $\exists \bar{z} \in g(\bar{x}) \cap (-Z_+)$ for which $\eta(\bar{z}) \leq 0$. Hence $\eta(\bar{z}) = 0$, which means that

$$0 \in \eta(g(\bar{x})). \qquad (3.14)$$

Know that $x \in D$ means $O \in h(x)$; and $\exists z \in g(x) \cap (-Z_+)$, which yields $\eta(z) \leq 0$. These and (3.13) deduce that $\xi(f(x) - \bar{y}) \geq 0, \forall x \in D$. According to the Slater constraint qualification, we have $\xi \neq O$. So we may take $y_0 \in \operatorname{int} Y_+$ such that $\xi(y_0) = 1$. Define the operator $S : Z \to Y$ and $T : W \to Y$ by

$$S(z) = \eta(z) y_0, T(w) = \varsigma(w) y_0 . \qquad (3.15)$$

It is easy to see that

$$S \in B^+(Z,Y), S(Z_+) = \eta(Z_+) y_0 \subseteq Y_+, T \in B(W,Y).$$

And (3.14) implies that

$$S(g(\bar{x})) = \eta(g(\bar{x}))y_0 \in 0 \cdot Y_+ = O. \quad (3.16)$$

Since $\bar{x} \in D$, we have $O \in h(\bar{x})$. Hence

$$O \in T(h(\bar{x})). \quad (3.17)$$

Therefore by (3.16) and (3.17)

$$\bar{y} \in f(\bar{x}) \subseteq f(\bar{x}) + S(g(\bar{x})) + T(h(\bar{x})).$$

From (3.13) and (3.14)

$$\begin{aligned} &\xi[f(x) + S(g(x)) + T(h(x))] \\ &= \xi(f(x)) + \eta((g(x))\xi(y_0) + \varsigma(h(x))\xi(y_0) \\ &= \xi(f(x)) + \eta(g(x)) + \varsigma(h(x)) \\ &\geq \xi(\bar{y}), \forall x \in D, \end{aligned}$$

i.e.,

$$\xi[f(x) - \bar{y}) + S(g(x)) + T(h(x))] \geq 0, \forall x \in D. \quad (3.18)$$

Take $F(x) = f(x) + S(g(x)) + T(h(x))$, $G(x) = \{O\}$ and $H(x) = \{O\}$, applying Theorem 1.2 to the functions $F(x) - \bar{y}, G(x), H(x)$, then (3.18) deduces that

$$(\bar{y} - [f(D) + S(g(D)) + T(h(D))] \cap \operatorname{int} Y_+ = \varnothing. \quad (3.19)$$

Hence

$$\bar{y} \in P\min[L(X, \bar{S}, \bar{T}), Y_+]. \quad (3.20)$$

On the other hand, since $O \in \bar{S}(g(\bar{x}))$ ((3.16)), $\exists \bar{z} \in g(\bar{x})$ for which $\bar{S}(\bar{z}) = O$. This and $(f(\bar{x}) - \bar{y}) \cap \operatorname{int} Y_+ = \varnothing$ together deduce that

$$(f(\bar{x}) - \bar{y} - S(\bar{z})) \cap \operatorname{int} Y_+ = \varnothing. \quad (3.21)$$

Combining the assumption $-g(\bar{x}) \subseteq Z_+$, $h(\bar{x}) = \{O\}$ and (3.20), (3.21), by Theorem 3.1 we conclude that $(\bar{x}, \bar{S}, \bar{T}) \in X \times B^+(Z, Y) \times B(W, Y)$ is a vector saddle-point of $L$. □

**Remark 3.1** For vector-valued functions, the condition $g(\bar{x}) \subseteq -Z_+$, $h(\bar{x}) = \{O\}$ in Theorem 3.2, 3.3 are always satisfied if $\bar{x} \in D$.

## 4. Scalar Saddle-Point Theorems

**Definition 4.1** Given $\bar{\xi} \in Y_+^* \setminus \{O\}$. The real-valued Lagrangian function of (VP) $l_{\bar{\xi}} : X \times Z_+^* \times W^* \to R$ is defined by

$$l_{\bar{\xi}}(x,\eta,\varsigma) = \bar{\xi}(f(x)) + \eta(g(x)) + \varsigma(h(x)).$$

**Definition 4.2** Given $\bar{\xi} \in Y_+^* \setminus \{O\}$. A triple $(\bar{x},\bar{\eta},\bar{\varsigma})$ is said to be a scalar saddle-point of the Lagrangian function $l_{\bar{\xi}}$, if

$$l_{\bar{\xi}}(\bar{x},\eta,\varsigma) \le l_{\bar{\xi}}(\bar{x},\bar{\eta},\bar{\varsigma}) \le l_{\bar{\xi}}(x,\bar{\eta},\bar{\varsigma}),$$

$\forall x \in D, \forall(\eta,\varsigma) \in Z_+^* \times W^*$.

The definition of a scalar saddle-point is the common definition of a saddle-point. To compare with the concept of a vector saddle-point in the previous section of this paper, we call it here a scalar saddle-point.

**Theorem 4.1** Suppose $\bar{x} \in D$, and

(a) $f$: $X \to 2^Y$ is $Y_+$ - preconvexlike on $D$, and $g$: $Y \to 2^Z$ is $Z_+$ - preconvexlike on $D$,
(b) $h$: $X \to 2^W$ is preaffine on $D$, and $int[\, h(D)] \neq \varnothing$,
(c) (VP) satisfies the Slater constrained qualification (SC).

If $\bar{x} \in D$ is a weakly efficient solution of (VP) for which $g(\bar{x}) \subseteq -Z_+$, $h(\bar{x}) = \{O\}$, then $\exists(\bar{\xi},\bar{\eta},\bar{\varsigma}) \in (Y_+^* \setminus \{O\}) \times Z_+^* \times W^*$ such that $(\bar{x},\bar{\eta},\bar{\varsigma})$ is a scalar saddle-point of the Lagrangian function $l_{\bar{\xi}}$ and $\bar{\eta}(g(\bar{x})) = \{0\}$.

**Proof**. Suppose that $\bar{x} \in D$ is a weakly efficient solution of (VP). Similar to the proof of (3.18) in Theorem 3.3 (in fact, $f(\bar{x})$ takes place of $\bar{y}$ from (3.12) to (3.18)) $\exists(\bar{S},\bar{T}) \in B^+(Z,Y) \times B(W,Y)$ such that

$$O \in \bar{S}(g(\bar{x})), \qquad (4.1)$$

and ( note that $g(\bar{x}) \subseteq -Z_+$, $h(\bar{x}) = \{O\}$ )

$$
\begin{aligned}
&\bar{\xi}(f(\bar{x}))+\bar{S}(g(\bar{x}))+\bar{T}(h(\bar{x}))) \\
&\leq \bar{\xi}(f(\bar{x})) \\
&\leq \bar{\xi}[f(x))+\bar{S}(g(x))+\bar{T}(h(x))] \\
&= \bar{\xi}(f(x))+\bar{\xi}\circ\bar{S}(g(x))+\bar{\xi}\circ\bar{T}(h(x)), \forall x\in D.
\end{aligned} \tag{4.2}
$$

Take $\bar{\eta}=\bar{\xi}\circ\bar{S}, \bar{\varsigma}=\bar{\xi}\circ\bar{T}$, then $(\bar{\eta},\bar{\varsigma})\in Z_+^*\times W^*$. And from (4.2) we attain

$$
\begin{aligned}
&l_{\bar{\xi}}(\bar{x},\bar{\eta},\bar{\varsigma}) \\
&\leq \bar{\xi}(f(\bar{x}))+\bar{\eta}(g(\bar{x}))+\bar{\varsigma}(h(\bar{x})) \\
&= \bar{\xi}[f(\bar{x})+\bar{S}(g(\bar{x}))+\bar{T}(h(\bar{x})] \\
&\leq \bar{\xi}[f(x)+\bar{S}(g(x))+\bar{T}(h(x)] \\
&= \bar{\xi}(f(x))+\bar{\eta}(g(x))+\bar{\varsigma}(h(x)) \\
&= l_{\bar{\xi}}(x,\bar{\eta},\bar{\varsigma}), \forall x\in D.
\end{aligned} \tag{4.3}
$$

And

$$
\bar{\xi}(f(\bar{x}))+\bar{\eta}(g(\bar{x}))\leq \bar{\xi}(f(x))+\bar{\eta}(g(x))+\bar{\varsigma}(h(x)), \forall x\in D. \tag{4.4}
$$

Since $x\in D$ we have $0\in\bar{\varsigma}(h(x))$, and (4.1) states that $0\in\bar{\eta}(g(\bar{x}))$, therefore from (4.4)

$$
\bar{\xi}(f(\bar{x}))\leq \bar{\xi}(f(x))+\bar{\eta}(g(x)), \forall x\in D.
$$

Taking $x=\bar{x}$ in above we attain $0\leq\bar{\eta}(g(\bar{x}))$. However, the assumption $g(\bar{x})\subseteq -Z_+$ implies that $\bar{\eta}(g(\bar{x}))\leq 0$. Thus

$$
\bar{\eta}(g(\bar{x}))=\{0\}. \tag{4.5}
$$

Hence from (4.5), noting that $\eta(g(\bar{x}))\leq 0$, and $h(\bar{x})=\{O\}$, we obtain

$$
\begin{aligned}
&l_{\bar{\xi}}(\bar{x},\eta,\varsigma) \\
&= \bar{\xi}(f(\bar{x}))+\eta(g(\bar{x}))+\varsigma(h(\bar{x})) \\
&\leq \bar{\xi}(f(\bar{x}))+\bar{\eta}(g(\bar{x}))+\bar{\varsigma}(h(\bar{x})) \\
&= l_{\bar{\xi}}(\bar{x},,\bar{\eta},\bar{\varsigma}), \forall(\eta,\varsigma)\in Y_+^*\times W^*.
\end{aligned} \tag{4.6}
$$

Combining (4.3) and (4.6), $(\bar{x},\bar{\eta},\bar{\varsigma})$ is a saddle-point of the Lagrangian function $l_{\bar{\xi}}$. □

**Remark 4.1** For vector-valued functions, the condition $g(\bar{x})\subseteq -Z_+$, $h(\bar{x})=\{O\}$ in Theorem 4.1 are always satisfied if $\bar{x}\in D$.

**Theorem 4.2** Let $\bar{x} \in D$. If $\exists(\bar{\xi},\bar{\eta},\bar{\varsigma}) \in (Y_+^* \setminus \{O\}) \times Z_+^* \times W^*$ for which $(\bar{x},\bar{\eta},\bar{\varsigma})$ is a saddle-point of the Lagrangian function $l_{\bar{\xi}}$, then $\bar{x} \in D$ is a weakly efficient solution of (VP) and $\bar{\eta}(g(\bar{x})) = \{0\}$, $\bar{\varsigma}(h(\bar{x})) = \{0\}$.

**Proof**. Suppose $\exists(\bar{\xi},\bar{\eta},\bar{\varsigma}) \in (Y_+^* \setminus \{O\}) \times Z_+^* \times W^*$ such that $(\bar{x},\bar{\eta},\bar{\varsigma})$ is a scalar saddle-point of the Lagrangian function $l_{\bar{\xi}}$, i.e.,

$$l_{\bar{\xi}}(\bar{x},\eta,\varsigma) \le l_{\bar{\xi}}(\bar{x},\bar{\eta},\bar{\varsigma}) \le l_{\bar{\xi}}(x,\bar{\eta},\bar{\varsigma}), \forall(\eta,\varsigma) \in Z_+^* \times W^*, \ \forall x \in D.$$

That is to say

$$\begin{aligned}&\bar{\xi}(f(\bar{x})) + \eta(g(\bar{x})) + \varsigma(h(\bar{x}))) \\ &\le \bar{\xi}(f(\bar{x})) + \bar{\eta}(g(\bar{x})) + \bar{\varsigma}(h(\bar{x}))), \forall x \in D. \ \forall(\eta,\varsigma) \in Z_+^* \times W^* \quad .\end{aligned}$$

Then

$$\eta(g(\bar{x})) + \varsigma(h(\bar{x}))) \le \bar{\eta}(g(\bar{x})) + \bar{\varsigma}(h(\bar{x})), \ \forall(\eta,\varsigma) \in Z_+^* \times W^*.$$

On the other hand

$$\bar{\xi}(f(\bar{x})) + \bar{\eta}(g(\bar{x})) + \bar{\varsigma}(h(\bar{x}))) \le \bar{\xi}(f(x)) + \bar{\eta}(g(x)) + \bar{\varsigma}(h(x)), \forall x \in D. \qquad (4.7)$$

Take $\eta = \bar{\eta}$, or $\varsigma = \bar{\varsigma}$ in (4.7) we have

$$\begin{aligned}&\eta(g(\bar{x})) \le \bar{\eta}(g(\bar{x})), \forall \eta \in Z_+^*, \\ &\varsigma(h(\bar{x}))) \le \bar{\varsigma}(h(\bar{x})), \forall \varsigma \in W^*.\end{aligned}$$

Therefore, taking $\eta = O$ in we get $\bar{\eta}(g(\bar{x})) \ge 0$, but taking $\eta = 2\bar{\eta}$ in we get $\bar{\eta}(g(\bar{x})) \le 0$. Hence

$$\bar{\eta}(g(\bar{x})) = \{0\}. \qquad (4.8)$$

Similarly,

$$\bar{\varsigma}(h(\bar{x})) = \{0\}. \qquad (4.9)$$

Noting that $g(x) \cap (-Z_+) \ne \varnothing$ and $0 \in \bar{\varsigma}(h(x))$ (since $x \in D$), according to (4.7), (4.8), and (4.9) we obtain

$$\bar{\xi}(f(\bar{x})) \le \bar{\xi}(f(x)), \forall x \in D.$$

Therefore, by Theorem 2.1, $\bar{x}$ is a weakly efficient solution of (VP). □

## References

[1] K. Fan, Minimax Theorem, Proceedings of the National Academy of Sciences of USA, 39(1953), 42-47.

[2] J. Borwein, Multivalued Convexity and Optimization: a Unified Approach to Inequality and equality Constraints, Mathematical Programming 13(1977), 183-199.

[3] Z. F. Li and S. Y. Wang, Connectedness of Super Efficient Sets in Vector Optimization of Set-Valued Maps, Mathematics of Operational Research. 48(1998), 207-217.

[4] V. Jeyakumar, A Generalization of a Minimax Theorem of Fan via a Theorem of the Alternative, Journal of Optimization Theory and Applications, 48(1986), 525-533.

[5] V. Jeyakumar, Convexlike Alternative Theorems and Mathematical Programming, Optimization, 16 (1985), 643- 652.

[6] Z. Li, The Optimality Conditions for Vector Optimization of Set-Valued Maps, Journal of Mathematical Analysis and Applications, 237(1999), 413-424.

[7] Z, Li, The Optimality Conditions of Differentiable Vector Optimization Problems, Journal of Mathematical Analysis and Applications, 201(1996), 35-43.

[8] V. Jeyakumar, A General Farkas Lemma and Characterization of Optimality for a Nonsmooth Program Involving Convex Processes, Journal of Optimization Theory and Application, 55(1987), 449-467.

[9] L. J. Lin, Optimization of Set-Valued Functions, Journal of Mathematical Analysis and Applications, 186(1994), 30-51.

[10] D. T. Luc, On Duality Theorem in Multiobjective Programming, Journal of Optimization Theory and Application, 48(1984): 557-582.

[11] X. M. Yang, X. Q. Yang and G. Y. Chen, Theorem of Alternative and Optimization, Journal of Optimization Theory and Application, 107(2000)3, 627-640.

[12] S. Paeck, Convexlike and Concavelike Conditions in Alternative, Minimax, and Minimization Theorems, Journal of Optimization Theory and Applications, 74(1992)2, 317-332.

[13] Y. Chen, G. Lan, and Y. Ouyang, Optimal Primal-Dual Methods for a Class of Saddle Point Problems, SIAM Journal on Optimization, 24(2014)4, 1779-1814.

[14] Q. B. Zhang, M. J. Li and C. Z. Cheng, Generalized Saddle Points Theorems for Set-Valued Mappings in Locally Generalized Convex Spaces, Nonlinear Analysis: Theory, Methods & Applications, 71(2009)1-2, 212-218.

[15] K. Q. Zhao and X. M. Yang, E-Propper Saddle Points and E-Proper Duality in Vector Optimization with Set-Valued Maps, Taiwanese Journal of Mathematics, 18(2014)2, 483-495.

[16] Y. Zhang and T. Chen, Minimax Problems for Set-Valued Mappings with Set Optimization, Numerical Algebra Control and Optimization, 4(2014)4, 327-340.

[17] P. Balamurugan, F. B. Inria and E. N. Supérieure, Stochastic Variance Reduction Methods for Saddle-Point Problems, 29th Conference on Neural Information Processing Systems (NIPS 2016), Barcelona, Spain.

[18] Y. Zhang and S. Li, Minimax Problems of Uniformly Same-Order Set-Valued Mappings, Bulletin of Korean Mathematical Society, 50(2013)5, 1639-1650.

[19] A. Kumar and P. K. Garg, Mixed Saddle Point and Its Equivalence with an Efficient Solution under Generalized (V, $p$)-Invexity, Applied Mathematics, 28(2015)9, 1630-1637.

[20] R. P. Agarwal, M. Balaj and D. O'Regan, An Intersection Theorem for Set-Valued Mappings, Applications of Mathematics, 58(2013)3, 269-278.

[21] G. Yu and X. Kong, Optimality and Duality in Set-Valued Optimization using Higher-order Radial Derivatives, Statistics Optimization and Information Computing, 4(2016) 154-162.

[22] C. Clason and T . Valkonen, Stability of Saddle Points via Explicit Coderivatives of Pointwise Subdifferentials, Set-Valued and Variational Analysis, 25(2017)1, 69-112.

[23] X. B. Li, Z. Lin, Q. L. Wang, and J. W. Chen, Holder Continuity of the Saddle Point Set for Real-Valued Functions, Numerical Functional Analysis and Optimization, 5(2017), 1-16.

[24] C. Wang, Q., Liu and B. Qu, Global Saddle Points of Nonlinear Augmented Lagrangian Functions, Journal of Global Optimization, 68(2017)1, 125-146.

[25] G. Yu, Strong Duality with Super Efficiency in Set-Valued Optimization, Journal of Nonlinear Science and Applications, 10 (2017), 3261-3272.